\documentclass[11pt,a4paper,twoside]{article}
\usepackage{amsthm,amsfonts,amsmath}

\newtheorem{theorem}{Theorem}%[section]
\newtheorem{definition}{Definition}%[section]
\newtheorem{example}{Example}%[section]
\newtheorem{lemma}{Lemma}%[section]
\newtheorem{proposition}{Proposition}%[section]
\newtheorem{remark}{Remark}%[section]

\author{Vladimir Rovenski\footnote{Department of Mathematics, University of Haifa, 3498838 Haifa, Israel.
       \newline e-mail: \texttt{vrovenski@univ.haifa.ac.il}\,.
       \quad ORCID: 0000-0003-0591-8307} }

\title{Foliated structure of weak nearly Sasakian manifolds}

\begin{document}

\date{}

\maketitle

\begin{abstract}
Weak contact metric manifolds, i.e., the linear complex structure on the contact distribu\-tion is replaced by a nonsingular skew-symmetric tensor,
defined by the author and R.\,Wolak, allowed a new look at the theory of contact manifolds.
In this paper we study the new structure of this type, called the weak nearly Sasakian structure.
We find conditions that are satisfied by almost contact manifolds and under which
the contact distribution is curvature invariant and the weak nearly Sasakian structure foliates into two types of totally geodesic foliations.
Our main result generalizes the theorem by B.\,Cappelletti-Montano and G.\,Dileo (2016) about the foliated structure of nearly Sasakian manifolds.


\end{abstract}

%\setcounter{secnumdepth}{4}
%%%%%%%%%%%%%%%%%%%%%%%%%%%%%%%%%%%%%%%%%%
%\setcounter{page}{1}

\section{Introduction}
\label{sec:00-ns}
%\label{sec1}

Riemannian contact geometry plays an important role in both modern mathematics and mechanics in explaining physical processes.
A popular class of almost contact metric manifolds $M^{\,2n+1}(\varphi,\xi,\eta,g)$ consists of Sasakian manifolds characterized by the equality, see \cite{blair2010riemannian},
\[
 (\nabla_X\,\varphi)Y=g(X,Y)\,\xi -\eta(Y)X,\quad X,Y\in\mathfrak{X}_M .
\]
Here, $g$ is a Riemannian metric, $\nabla$ is the Levi-Civita connection,
$\varphi$ is a $(1,1)$-tensor, $\xi$ is a Reeb vector field and $\eta$ is a 1-form, satisfying $\eta(\xi)=1$ and
\[
 g(\varphi X,\varphi Y)= g(X,Y) -\eta(X)\,\eta(Y),\quad X,Y\in\mathfrak{X}_M ,
 %\varphi^2 = -{\rm id}_{\,TM} + \eta\otimes \xi, \qquad .
\]
and $\mathfrak{X}_M$ is the space of vector fields on $M$.
D.\,Blair, D.\,Showers and Y.\,Komatu \cite{blair1976} defined nearly Sasakian structure $(\varphi,\xi,\eta,g)$
using a similar condition for the symmetric part of $\nabla\varphi$:
\begin{equation}\label{E-nS-02}
 (\nabla_Y\,\varphi)Y = g(Y,Y)\,\xi -\eta(Y)Y,\quad Y\in\mathfrak{X}_M,
\end{equation}
and showed that a normal nearly Sasakian structure is Sasakian and hence is contact.
Any 3-dimensional nearly Sasakian manifold is Sasakian, see \cite{Ol-1980},
and any 5-dimensional nearly Sasakian manifold has Einstein metric of positive scalar curvature, see~\cite{C-MD-2016};
an example is a sphere $S^{\,5}$ with the almost contact metric structure induced by the almost Hermitian structure of $S^{\,6}$.
%(defined by the cross product of the imaginary part of the octonions).

The Reeb vector field $\xi$ of a nearly Sasakian structure is a unit Killing vector field (an infinitesimal generator of isometries or symmetries).
The~influence of constant-length Killing vector fields on the Riemannian geometry has been studied by many authors, e.g.,~\cite{N-2021}.

B.\,Cappelletti-Montano and G.\,Dileo proved in \cite{C-MD-2016} that nearly Sasakian manifolds foliate into two types of
totally geodesic foliations.
Totally geodesic foliations (e.g., fibrations or submersions with totally geodesic fibers)
of Riemannian manifolds have the simplest extrinsic geometry of the leaves and appear
%naturally
as kernels of degenerate tensors, e.g., \cite{rov-2023c,Rov-Wa-2021}.
A.D.\,Nicola, G.\,Dileo and I.\,Yudin in \cite{NDY-2018}, using results of \cite{C-MD-2016,Ol-1979},
gave a criterion for the almost contact metric structure to be Sasakian; namely,
that every nearly Sasakian manifold of dimension greater than 5 is~Sasakian.

In \cite{RP-2,RWo-2,rov-117}
and \cite[Section~5.3.8]{Rov-Wa-2021}, we introduced and studies metric structures on a smooth manifold that genera\-lize the almost contact, Sasakian, etc. metric structures.
These so-called ``weak" structures (the~li\-near complex structure on the contact distribution is replaced by a nonsingular skew-symmetric tensor)
made it possible to take a new look at the classical structures and find new applications.
In~\cite{rov-2023} we investigated new structures of this type, called {weak nearly Sasakian structure} and {weak nearly cosymplectic structure}.

In this~paper we continue our study of the weak nearly Sasakian structure.
We find conditions \eqref{E-nS-10} and \eqref{E-nS-04c} that are satisfied by almost contact manifolds and under which
the Reeb vector field is Killing, the contact distribution is curvature invariant and the weak nearly Sasakian structure foliates into two types of
totally geodesic foliations.
Section~\ref{sec:01-ns}, following the introductory Section~\ref{sec:00-ns}, reviews the basics of weak almost contact manifolds.
Section~\ref{sec:02-ns} contains auxiliary results on the geometry of weak almost Sasakian structure;
in particular, we find the expression of the Ricci tensor in the $\xi$-direction
and show that the contact distribution $\ker\eta$ is curvature invariant under conditions \eqref{E-nS-10} and \eqref{E-nS-04c}.
In~Section~\ref{sec:04-ns}, we prove that a weak almost Sasakian manifold satisfying \eqref{E-nS-10} and \eqref{E-nS-04c}, is foliated.
Our main result (Theorem~\ref{T-2.2}) generalizes Theorems~3.3 and 3.5 of \cite{C-MD-2016}.
In~Section~\ref{sec:app}, using the approach of \cite{rov-2023b} we prove an auxiliary proposition necessary for the main result.
Our~proofs use the properties of new tensors, as well as classical~constructions.

\section{Preliminaries}
\label{sec:01-ns}

A~\textit{weak almost contact structure} on a smooth manifold $M^{\,2n+1}\ (n\ge1)$ is a set $(\varphi,Q,\xi,\eta)$,
where $\varphi$ is a $(1,1)$-tensor, $\xi$ is a vector field (called Reeb vector field), $\eta$ is a
1-form and $Q$ is a nonsingular $(1,1)$-tensor on $TM$, satisfying, see \cite{RP-2,RWo-2},
\begin{equation}\label{E-nS-2.1}
 \varphi^2 = -Q + \eta\otimes \xi, \quad \eta(\xi)=1,\quad Q\,\xi=\xi .
\end{equation}
By \eqref{E-nS-2.1}, $\eta$ defines a smooth $2n$-dimensional distribution $\ker\eta$.
We assume that $\ker\eta$ is $\varphi$-invariant, i.e., $\varphi(\ker\eta)\subset\ker\eta$
(as in the classical theory \cite{blair2010riemannian}, where $Q={\rm id}_{\,TM}$).
By this and \eqref{E-nS-2.1}, $\ker\eta$ is $Q$-invariant, i.e., $Q(\ker\eta)\subset\ker\eta$,
and the following is true:
\[
 \varphi\,\xi=0,\quad \eta\circ\varphi=0,\quad \eta\circ Q=\eta,\quad [Q,\,\varphi]:={Q}\circ\varphi - \varphi\circ{Q}=0.
\]
A ``small" (1,1)-tensor $\widetilde Q= Q - {\rm id}_{\,TM}$ is a measure
of the difference between a weak almost contact structure and an almost contact one.
Note that
\[
 [\widetilde{Q},\varphi]:=\widetilde{Q}\circ\varphi - \varphi\circ\widetilde{Q}=0,\quad
 \eta\circ\widetilde Q=0,\quad
 \widetilde{Q}\,\xi=0.
\]
If there is a Riemannian metric $g$ on $M$ such that
\begin{align}\label{E-nS-2.2}
g(\varphi X,\varphi Y)= g(X,Q\,Y) -\eta(X)\,\eta(Y),\quad X,Y\in\mathfrak{X}_M,
\end{align}
then $(\varphi,Q,\xi,\eta,g)$ is called a {\it weak almost contact metric structure} on $M$.
A weak almost contact manifold $M^{\,2n+1}(\varphi,Q,\xi,\eta)$ endowed with a compatible Riemannian metric $g$ is called a \textit{weak almost contact metric manifold} and is denoted by $M^{\,2n+1}(\varphi,Q,\xi,\eta,g)$.

For a weak almost contact metric structure, $\varphi$ is skew-symmetric and $Q$ is self-adjoint.
Taking covariant derivative of the equality $g(\varphi Y,Z)=-g(Y,\varphi Z)$, we see that $\nabla_{X}\,\varphi$ is skew-symmetric:
\begin{equation*}%\label{E-nS-05e}
 g((\nabla_{X}\,\varphi) Y, Z)=-g((\nabla_{X}\,\varphi) Z, Y).
\end{equation*}
Setting $Y=\xi$ in \eqref{E-nS-2.2}, we get $\eta(X)=g(\xi, X)$, as in the classical theory.
By \eqref{E-nS-2.2},
\[
 g(X,Q\,X)=g(\varphi X,\varphi X)>0
\]
is true for any nonzero vector $X\in\ker\eta$; thus, the tensor $Q$ is positive definite.

A~\textit{weak contact metric structure} is a weak almost contact metric structure satisfying
\[
 d\eta(X,Y)=g(X,\varphi Y),\quad X,Y\in\mathfrak{X}_M,
\]
where the exterior derivative $d\eta$ of $\eta$ is given~by (for example,~\cite{blair2010riemannian})
\begin{align*}
 d\eta(X,Y) = \frac12\,\{X(\eta(Y)) - Y(\eta(X)) - \eta([X,Y])\}.
\end{align*}
A~1-form $\eta$ on a manifold $M^{\,2n+1}$ is said to be \textit{contact} if $\eta\wedge (d\eta)^n\ne0$, e.g.,~\cite{blair2010riemannian}.
Recall that for a weak contact metric structure $(\varphi,Q,\xi,\eta,g)$, the 1-form $\eta$ is contact, see \cite[Lemma~2.1]{rov-2023b}.

\begin{definition}[see \cite{rov-2023}]\rm
A weak almost contact metric manifold $M^{\,2n+1}(\varphi,Q,\xi,\eta,g)$ is called \textit{weak nearly Sasaki\-an}~if
\eqref{E-nS-02} is true, or, equivalently,
\begin{equation}\label{E-nS-00b}
 (\nabla_Y\,\varphi)Z + (\nabla_Z\,\varphi)Y = 2\,g(Y,Z)\,\xi -\eta(Z)Y -\eta(Y)Z,\qquad Y,Z\in\mathfrak{X}_M.
\end{equation}
\end{definition}

In addition to \eqref{E-nS-00b}, the following two conditions for weak almost contact manifolds that are automatically satis\-fied by almost contact manifolds ($\widetilde Q=0$)
play the key role in this paper:
\begin{align}\label{E-nS-10}
 & (\nabla_X\,\widetilde Q)\,Y=0,\quad X\in\mathfrak{X}_M,\ Y\in\ker\eta, \\
\label{E-nS-04c}
 & R_{\widetilde Q X,Y}Z\in\ker\eta,\quad X,Y,Z\in\ker\eta ,
\end{align}
where the curvature tensor $R$ is given by, e.g., \cite{CLN-2006},
\[
 R_{{X},{Y}}Z=\nabla_X\nabla_Y Z -\nabla_Y\nabla_X Z -\nabla_{[X,Y]} Z,\quad X,Y,Z\in\mathfrak{X}_M.
\]
By \eqref{E-nS-04c} and the first Bianchi identity, e.g., \cite{CLN-2006}:
%\[
 $R_{\,X,Y}\,\widetilde Q Z\in\ker\eta$ for all $X,Y,Z\in\ker\eta$.
%\]

\begin{example}\rm
Let $M(\varphi,Q,\xi,\eta,g)$ be a three-dimensional weak nearly Sasakian ma\-nifold with the~condition \eqref{E-nS-10}.
By \eqref{E-nS-2.1}, the symmetric tensor $Q$ has on the plane field $\ker\eta$ the form $\lambda\,{\rm id}_{\,\ker\eta}$ for some positive $\lambda\in\mathbb{R}$.
It was shown in \cite{rov-2023} that this structure is reduced to the nearly Sasakian structure $(\tilde\varphi,\xi,\eta,\tilde g)$ on $M$, where
\begin{align*}
 \tilde\varphi = \lambda^{\,-\frac12}\,\varphi ,\quad
 \tilde  g|_{\,\ker\eta} = \lambda^{\,\frac12}\,g|_{\,\ker\eta},\quad
 \tilde g(\xi,\,\cdot) = {g}(\xi,\,\cdot) .
\end{align*}
Since $\dim M=3$, the
%nearly Sasakian
structure $(\tilde\varphi,\xi,\eta,\tilde g)$ is Sasakian, see~\cite[Theorem~5.1]{Ol-1980}.
\end{example}

Using the equalities $(\nabla_\xi\,\varphi)\,\xi = 0$ and $\varphi\,\xi=0$ for a weak nearly Sasakian manifold, we conclude that $\xi$ is a geodesic vector field ($\nabla_\xi\,\xi=0$).
Recall \cite{rov-2023}~that on a weak nearly Sasakian manifold with the property \eqref{E-nS-10} the unit vector field $\xi$ is Killing ($\pounds_\xi\,g=0$, e.g., \cite{CLN-2006}). Here $\pounds$ is the Lie derivative and the following identity is true:
\[
 (\pounds_\xi\,g)(X,Y)=g(\nabla_X\,\xi, Y) + g(\nabla_Y\,\xi, X),\quad X,Y\in\mathfrak{X}_M.
\]
Therefore, $\xi$-curves determine a Riemannian geodesic foliation.
Note that $-\nabla\xi$ is its splitting operator (or, co-nullity tensor), see \cite[Section~1.3.1]{Rov-Wa-2021}.
Moreover, $\eta$ is a Killing 1-form:
%\[
 $(\nabla_X\,\eta)(X)=0$ for all $X\in\mathfrak{X}_M$.
%\]
Using \eqref{E-nS-10} and $\nabla_\xi\,\xi=0$, we obtain
\begin{equation}\label{E-nabla-xi-Q}
 \nabla_\xi\,\widetilde Q=0.
\end{equation}
For a Riemannian manifold equipped with a Killing vector field ${\xi}$, we get, see \cite{YK-1985},
\begin{equation}\label{E-nS-04}
 \nabla_X\nabla_Y\,{\xi} - \nabla_{\nabla_X Y}\,{\xi} = R_{\,X,\,{\xi}}\,Y .
\end{equation}

\section{Auxiliary results}
\label{sec:02-ns}

In this section we generalize some properties of almost Sasakian manifolds to the case of weak almost Sasakian manifolds $M^{\,2n+1}(\varphi,Q,\xi,\eta,g)$ with the conditions \eqref{E-nS-10} and~\eqref{E-nS-04c}.
Define a (1,1)-tensor field $h$ on $M$, as in the classical case, e.g., \cite{C-MD-2016},
\begin{equation}\label{E-c-01}
 h = \nabla\xi + \varphi.
\end{equation}
We get $\eta\circ h = 0$ and $h(\ker\eta)\subset\ker\eta$. Since $\xi$ is a geodesic vector field,
%($\nabla_\xi\,\xi=0$),
 we also get
\begin{equation}\label{E-cc-01}
 h\,\xi=0 .
\end{equation}
Since $\xi$ is a Killing vector field and $\varphi$ is skew-symmetric, the tensor $h$ is skew-symmetric:
\[
 g(h X,\, X) = g(\nabla_X\,\xi, X) + g(\varphi X, X) = \frac12\,(\pounds_\xi\,g)(X,X) = 0,
\]
and $\nabla_X\,\eta=g((h-\varphi)X,\,\cdot)$ is true.
The distribution $\ker\eta$ is integrable, $[X,Y]\in\ker\eta\ (X,Y\in\ker\eta)$, if and only if $h=\varphi$:
\[
 g([X,Y], \xi) = 2\,g((h-\varphi)Y, X),\quad X,Y\in\ker\eta;
\]
and in this case, our manifold splits along $\xi$ and $\ker\eta$ (is locally the metric product).

\begin{lemma}
%\label{L-nS-02}
For a weak nearly Sasakian manifold $M^{\,2n+1}(\varphi,Q,\xi,\eta,g)$ we obtain
\begin{align}
\label{E-nS-01a}
 & (\nabla_X\,h)\,\xi = -h(h-\varphi)X ,\\
\label{E-nS-01c}
 & (\nabla_X\,\varphi)\,\xi =  -\varphi(h-\varphi) X ,\\
%%%
\label{E-nS-01b}
 & h\,\varphi + \varphi\,h = -2\,\widetilde Q ,\\
\label{E-nS-01d}
 & h\,Q = Q\,h \quad (h\ {\rm commutes\ with}\ Q).
\end{align}
Moreover,
\begin{align}\label{E-nS-01e}
 h^2\varphi=\varphi h^2,\quad h\varphi^2=\varphi^2 h,\quad h^2\varphi^2=\varphi^2 h^2.
\end{align}
\end{lemma}

\begin{proof}
Differentiating the equality $h\,\xi=0$ and using the definition \eqref{E-c-01}, we get \eqref{E-nS-01a}:
\[
 0=\nabla_X\,(h\,\xi) = (\nabla_X\,h)\,\xi +h(\nabla_X\,\xi) =(\nabla_X\,h)\,\xi +h(h-\varphi) X.
\]
Using $\varphi\,\xi=0$ and the definition \eqref{E-c-01}, we get \eqref{E-nS-01c}:
\[
 (\nabla_X\,\varphi)\,\xi = -\varphi(\nabla_X\,\xi)= -\varphi(h-\varphi)X.
\]
Differentiating the equality $g(\varphi\,Y,\xi)=0$ yields
\[
 0=X g(\varphi\,Y,\xi) = g((\nabla_X\,\varphi)Y, \xi) +g(\varphi\,Y, (h-\varphi)X).
\]
Summing this with the equality
%\[
 $g((\nabla_Y\,\varphi)X, \xi) +g(\varphi\,X, (h-\varphi)Y)=0$
%\]
and applying \eqref{E-nS-00b}, gives~\eqref{E-nS-01b}.
By \eqref{E-nS-01b} and \eqref{E-nS-2.1}, using the equalities \eqref{E-cc-01}, $\eta\circ h=0$ and $[\widetilde Q, \varphi]=0$, we get \eqref{E-nS-01d}.
Using \eqref{E-nS-01b} and $[Q,\,\varphi]=0$, we obtain \eqref{E-nS-01e}.
\end{proof}

%\begin{remark}\rm
The contact distribution of a nearly Sasakian manifold is {curvature invariant}:
\begin{equation}\label{E-nS-04cc}
 R_{\,X,Y}Z\in\ker\eta,\quad X,Y,Z\in\ker\eta ,
\end{equation}
see~\cite{Ol-1979}.
For example, the tangent bundle of a totally geodesic submanifold in a Riemannian manifold and the distribution $\ker\eta$ of any 1-form $\eta$ on a real space form satisfy \eqref{E-nS-04cc}.
%For example, a real space form with a 1-form $\eta$ and the tangent bundle of any totally geodesic submanifold satisfy \eqref{E-nS-04cc}.
%
According to the following proposition, a weak nearly Sasakian manifold satisfies \eqref{E-nS-04cc}
if we assume a weaker condition~\eqref{E-nS-04c}.
%We will prove, see \eqref{E-nS-04ccc} of Proposition~\ref{P-R01}, that a weak nearly Sasakian manifold satisfies \eqref{E-nS-04cc} if we assume a weaker %condition~\eqref{E-nS-04c}.
%\end{remark}

\begin{proposition}\label{P-R01}
If a weak nearly Sasakian manifold
%$M(\varphi,Q,\xi,\eta,g)$
satisfies the condition \eqref{E-nS-04c}, then
\begin{align}\label{E-nS-04ccc}
 & g(R_{\,\xi, Z}\,\varphi X,\varphi Y) = 0,\quad {\rm hence,}\ \ker\eta\ {\rm is\ a\ curvature\ invariant\ distribution}.
\end{align}
\end{proposition}

The {Ricci tensor} ${\rm Ric}$ on $(M, g)$ is defined as the suitable trace of the curvature tensor:
\begin{equation*}
%\label{eq:ricci}
 {\rm Ric}\,(X,Y) = {\rm trace}_{\,g}(Z\to R_{\,Z,X}\,Y) = \sum\nolimits_{\,i} g(R_{\,E_i,X}\,Y, E_i)
\end{equation*}
when $(E_i)$ is an orthonormal frame.
%
% The following lemma is similar to Lemma~3.5 in~\cite{rov-2023b}.

\begin{lemma}\label{L-nS-04}
For a weak nearly Sasakian manifold $M^{\,2n+1}(\varphi,Q,\xi,\eta,g)$ with the condition \eqref{E-nS-04c}, we obtain
\begin{align}
\label{E-3.23}
 & R_{\,\xi, X}Y = -(\nabla_X\,(h-\varphi)\,)\,Y, \\
%%%%%%
\label{E-3.24}
 & (\nabla_X\,(h-\varphi))Y =  g( (h-\varphi)^2 X, Y)\,\xi - \eta(Y)\,(h-\varphi)^2 X , \\
%%%%%%
\label{E-3.25}
 & {\rm Ric}\,(\xi, Z) = -\eta(Z)\,({\rm tr}\,(h^2 + \widetilde Q) - 2\,n).
\end{align}
In particular, ${\rm tr}(h^2+ \widetilde Q) = const$, ${\rm Ric}\,(\xi, \xi) = const\ge0$ and
\begin{align}\label{E-3.24b}
 \nabla_\xi\,h = \nabla_\xi\,\varphi = \varphi h + \widetilde Q.
\end{align}
\end{lemma}

\begin{proof}
From \eqref{E-nS-04} (since $\xi$ is a Killing vector field) and \eqref{E-c-01}, we find \eqref{E-3.23}.
Replacing $Y$ by $\varphi Y$ and $Z$ by $\varphi Z$ in
\[
 g(R_{\,\xi, X}Y, Z) = g((\nabla_X\,(\varphi-h))Y, Z),
\]
see \eqref{E-3.23}, and using \eqref{E-nS-04ccc}, we get $g((\nabla_X\,(h-\varphi))\,\varphi Y, \varphi Z) = 0$, i.e.,
\begin{equation}\label{E-nS-03f}
 g((\nabla_X\,(h-\varphi)) Y, Z) = 0\quad (Y, Z\in\ker\eta).
\end{equation}
Then, using \eqref{E-nS-03f}, we find the $\xi$-component of $(\nabla_X\,(h-\varphi))Y$:
\begin{align*}
 & g((\nabla_X\,(h-\varphi))Y, \xi) = g(\nabla_X(h\,Y-\varphi Y), \,\xi) \\
 & = -g((h-\varphi)Y, \,\nabla_X\,\xi) = -g((h-\varphi)X, (h-\varphi)Y) = g((h-\varphi)^2 X, Y),
\end{align*}
and $\ker\eta$-component of $(\nabla_X\,(h-\varphi))Y$:
\begin{align*}
 & g((\nabla_X\,(h-\varphi))Y, Z) = \eta(Y)\,g((\nabla_X\,(h-\varphi))\,\xi, \,Z) \\
 & = \eta(Y)\,g(h(\varphi-h) X, Z) +\eta(Y)\,g(\varphi(h-\varphi)X, Z)  \\
 & = -\eta(Y)\,g( (h-\varphi)^2 X, Z)
 \quad (Z\in\ker\eta),
\end{align*}
from which \eqref{E-3.24} follows.
From \eqref{E-3.24} with $X=\xi$ we find $\nabla_\xi\,(h-\varphi)=0$.

Let $\{e_i\}\ (i=1,\ldots,2n+1)$ be a local orthonormal frame on $M$ with $e_{\,2n+1}=\xi$.
Putting $X=Y=e_i$ in \eqref{E-3.24}, then using \eqref{E-nS-04cc} and summing over $i=1,\ldots,2n+1$, we get
\[
 {\rm Ric}\,(\xi, Z) = -\eta(Z)\,{\rm tr}\,((h-\varphi)^2)
\]
and ${\rm Ric}\,(\xi, \xi)\ge0$.
From the above and $(h-\varphi)^2 = h^2+\varphi^2 +2\,\widetilde Q$, we get~\eqref{E-3.25}.

Using ${\rm tr}(\varphi^2)=-{\rm tr}\,(\widetilde Q+2\,n)$, we get $X({\rm tr}(\varphi^2)) = - X({\rm tr}\,\widetilde Q)$.
Similarly,
\begin{align*}
 & X({\rm tr}(h^2)) = {\rm tr}(\nabla_X\,(h^2)) ={\rm tr}(h\nabla_X\,h + (\nabla_X\,h)h) \\
 & = {\rm tr}(h\nabla_X(h-\varphi) + (\nabla_X(h-\varphi))h) + {\rm tr}(h\nabla_X\,\varphi + (\nabla_X\,\varphi)h) ,
\end{align*}
where, in view of \eqref{E-nS-01b},
\begin{align*}
 {\rm tr}(h\nabla_X\,\varphi + (\nabla_X\,\varphi)h) & = {\rm tr}(\nabla_X(h\varphi +\varphi h))
 - 2\,{\rm tr}\big(\varphi\nabla_X(h-\varphi)\big)
% - {\rm tr}\big(\varphi\nabla_X(h-\varphi) + (\nabla_X(h{-}\varphi))\varphi\big)
 {-} X({\rm tr}(\varphi^2)) \\
 & = - X({\rm tr}\,\widetilde Q) - 2\,{\rm tr}\big(\varphi\nabla_X(h{-}\varphi)
% - {\rm tr}\big(\varphi\nabla_X(h{-}\varphi) + (\nabla_X(h{-}\varphi))\varphi
\big) .
\end{align*}
By \eqref{E-3.24} we get
\[
 {\rm tr}\big(h\nabla_X(h-\varphi) + (\nabla_X(h-\varphi))h\big)=0,\quad
 {\rm tr}\big(\varphi\nabla_X(h-\varphi)
 %+ (\nabla_X(h-\varphi))\varphi
 \big) =0.
\]
Therefore, $X({\rm tr}(h^2)) = - X({\rm tr}\,\widetilde Q)$ for all  $X\in\mathfrak{X}_M$.
This implies ${\rm tr}(h^2+ \widetilde Q) = const$.
By this and \eqref{E-3.23}, we also get ${\rm Ric}\,(\xi, \xi) = -{\rm tr}\,(h^2 + \widetilde Q) + 2\,n = const$.
By \eqref{E-3.24} with $X=\xi$ we get $\nabla_\xi\,h = \nabla_\xi\,\varphi$.
Then applying \eqref{E-nS-00b} to $(\nabla_\xi\,\varphi)X\ (X\in\mathfrak{X}_M)$, we complete the proof of \eqref{E-3.24b}.
\end{proof}

\begin{remark}\rm
Note that $(h-\varphi)^2 X = h^2 X +\varphi^2 X + 2\,\widetilde Q X$.
By \eqref{E-3.23}--\eqref{E-3.24}, we get
\begin{align}\label{E-3.23b}
%\nonumber
% & R_{\,\xi, X}Y =  \eta(Y)\,(h-\varphi)^2 X - g( (h-\varphi)^2 X , Y)\,\xi, \\
 g(R_{\,\xi, X}Y, Z) =  \eta(Y)\,g( (h-\varphi)^2 X , Z) - \eta(Z)\,g( (h-\varphi)^2 X, Y).
\end{align}
\end{remark}

\section{Main results}
\label{sec:04-ns}

Here, we prove that weak nearly Sasakian manifolds with conditions \eqref{E-nS-10} and~\eqref{E-nS-04c} have a foliated structure.
First, we will generalize Lemma 2.1 in \cite{Ol-1980},
which characterizes Sasakian manifolds among nearly Sasakian manifolds by the condition $h=0$.

\begin{proposition}
%\label{Th-4.1}
For a weak nearly Sasakian manifold with the property \eqref{E-nS-10}, the equality $h=0$ holds if and only if the manifold is Sasakian.
\end{proposition}

\begin{proof}
Let $h=0$, i.e., $\varphi = -\nabla\xi$, see \eqref{E-c-01}.
For every vector fields $X, Y$ in $\ker\eta$ we have
\begin{equation*}
%\label{E-c-01b}
 2\,d\eta(X, Y ) = g(\nabla_X\,\xi, Y) - g(\nabla_Y\,\xi, X) = 2\,g(X, \varphi Y),
\end{equation*}
thus, our manifold is contact.
From \eqref{E-nS-01b} we find $\widetilde Q=0$, thus our manifold is nearly Sasakian.
Finally, by Lemma~\ref{L-nS-04} and the Bianchi identity, we obtain
\[
 R_{\,X, Y}\,\xi = R_{\,\xi, Y}X - R_{\,\xi, X}Y = \eta(X)\,\varphi^2 Y - \eta(Y)\,\varphi^2 X
  = \eta(Y)\,X - \eta(X)\,Y
% $g(R_{\,\xi, X}Y, \xi) = g( X, Y)$
\]
for all $X,Y\in\ker\eta$.
Thus,
by \cite[Proposition~7.6]{blair2010riemannian},
%the $\xi$-sectional curvature is 1, and
our manifold is Sasakian.
%Finally, by \eqref{E-3.23b}, we obtain
% $g(R_{\,\xi, X}Y, \xi) = g( X, Y)$ for $X,Y\in\ker\eta$.
%Thus, the $\xi$-sectional curvature is 1, and our manifold is Sasakian.
\end{proof}

The following proposition generalizes \cite[Proposition~3.2]{C-MD-2016}.

\begin{proposition}\label{prop-4.1}
For a weak nearly Sasakian manifold with conditions \eqref{E-nS-10} and \eqref{E-nS-04c},
the eigenvalues (and their multiplicities) of the symmetric operator $h^2$
%, $\varphi^2$ and $(h-\varphi)^2$
are constant.
\end{proposition}

\begin{proof}
From \eqref{E-3.23b}, \eqref{E-nS-2.1} and Lemma~\ref{L-nS-04} we obtain the following generalization of equation (10) in~\cite{C-MD-2016}:
\begin{align}\label{E-nS-11}
\nonumber
 &\quad (\nabla_X\, h^2)Y = h(\nabla_X\, h)Y + (\nabla_X\, h) h Y \\
\nonumber
 & = h(\nabla_X\,\varphi)Y +(\nabla_X\,\varphi) h Y -\eta(Y)h\,(h-\varphi)^2 X +g((h-\varphi)^2 X, h Y)\,\xi \\
\nonumber
 & = -\varphi(\nabla_X\,h)Y -(\nabla_X\,h)\,\varphi Y -2(\nabla_X\,\widetilde Q)Y -\eta(Y)h\,(h-\varphi)^2 X  + g((h-\varphi)^2 X, hY)\,\xi \\
\nonumber
 & = -\varphi(\nabla_X\,\varphi)Y + \eta(Y)\varphi(h-\varphi)^2 X - (\nabla_X\,\varphi)\,\varphi Y -g((h-\varphi)^2 X, \varphi Y)\,\xi \\
\nonumber
 &\quad - 2(\nabla_X\,\widetilde Q)Y -\eta(Y)h\,(h-\varphi)^2 X +g((h-\varphi)^2 X, h Y)\,\xi \\
\nonumber
 & = (\nabla_X\,\widetilde Q)Y - g((h-\varphi)X, Y)\,\xi - \eta(Y)(h-\varphi)X + \eta(Y)\varphi(h-\varphi)^2 X  \\
\nonumber
 &\quad - g((h-\varphi)^2 X, \varphi Y)\,\xi - 2\,(\nabla_X\,\widetilde Q)Y + g((h-\varphi)^2 X, h Y)\,\xi  -\eta(Y)\,h (h-\varphi)^2 X \\
%%%%%%%%%%%%%
 & = - g\big( (h-\varphi)(h^2 + \widetilde Q) X, Y \big)\,\xi - \eta(Y)(h-\varphi) (h^2 + \widetilde Q) X - (\nabla_X\,\widetilde Q)Y.
\end{align}
Consider an eigenvalue $\mu$ of $h^2$ and a local unit vector field $Y\bot\,\xi$ such that $h^2Y = \mu Y$.
Applying \eqref{E-nS-11} for any nonzero vector fields $X\in\mathfrak{X}_M$ and using \eqref{E-nS-10}, we find
$g((\nabla_X\, h^2)Y, Y) =0 $, thus
\begin{align*}
 & 0 = g((\nabla_X\, h^2)Y, Y) = g(\nabla_X\,(h^2Y), Y) - g(h^2(\nabla_X\,Y), Y ) \\
 & \ \ = X(\mu)\,g(Y, Y) + \mu\,g(\nabla_X\,Y, Y) - g(\nabla_X\,Y, h^2Y) = X(\mu)\,g(Y, Y),
\end{align*}
which implies that $X(\mu) = 0$ for all $X\in\mathfrak{X}_M$.
\end{proof}

Since $h$ is skew-symmetric, the nonzero eigenvalues of $h^2$ are negative.
By Proposition~\ref{prop-4.1}, the spectrum of the self-adjoint operator $h^2$ has the~form
\begin{equation}\label{E-nS-11b}
  Spec(h^2) = \{0, -\lambda_1^2,\ldots -\lambda_r^2\} ,
\end{equation}
where $\lambda_i$ is a positive real number and $\lambda_i\ne \lambda_j$ for $i\ne j$.

In particular, ${\rm tr}\,(h^2) = const\le0$, and by Lemma~\ref{L-nS-04}, ${\rm tr}\,Q = const>0$.

Denote by $[\xi]$ the 1-dimensional distribution generated by $\xi$,
and by $D_0$ a smooth distribution of the eigenvectors of $h^2$ with eigenvalue 0 orthogonal to $\xi$.
Denote by $D_i$ a smooth distribution of the eigenvectors of $h^2$ with eigenvalue $-\lambda^2_i$.
Note that the distributions $D_0$ and $D_i\ (i = 1,\ldots, r)$ belong to $\ker\eta$ and are $\varphi$-invariant and $h$-invariant.
In particular, the eigenvalue $0$ has multiplicity $2p+1$ for some integer $p\ge0$.
If $X$ is a unit eigenvector of $h^2$ with eigenvalue $-\lambda^2_i$, then by \eqref{E-nS-01b} and \eqref{E-nS-01d},
$X, \varphi X, hX$ and $h\,\varphi X$ are nonzero eigenvectors of $h^2$ with eigenvalue $-\lambda^2_i$.

\begin{lemma}\label{L-04}
The tensors $h$ and $\widetilde Q$ of a weak nearly Sasakian manifold are equal to zero on $D_0$.
\end{lemma}

\begin{proof}
Let $h^2 X=0$ for some $X\ne0$. Since $g(hX,hX)= -g(h^2X,X)=0$, we get $h X=0$.
Similarly, using \eqref{E-nS-01e}, we get
\[
 g(h\varphi X,h\varphi X)= -g(h^2\varphi X,\varphi X)= -g(\varphi h^2 X,\varphi X)=0;
\]
hence, $h\varphi X=0$. Taking into account the above and \eqref{E-nS-01b}, we obtain $\widetilde Q X = 0$.
\end{proof}

The following theorem generalizes \cite[Theorems~3.3 and 3.5]{C-MD-2016}.

\begin{theorem}\label{T-2.2}
Let $M^{\,2n+1}(\varphi,Q,\xi,\eta,g)$ be a weak nearly Sasakian
%(non-Sasakian)
manifold with conditions \eqref{E-nS-10} and \eqref{E-nS-04c},
and let the spectrum of the self-adjoint operator $h^2$ have the form \eqref{E-nS-11b},
where the eigenvalue $0$ has multiplicity $2p+1$ for some integer $p\ge0$.
Then, the distribution $[\xi]\oplus D_0$ and
each distribution $[\xi]\oplus D_i\ (i = 1,\ldots, r)$ are integrable with totally geodesic leaves.
If $p>0$, then

\noindent\ \
$(a)$ the distribution $[\xi]\oplus D_1\oplus\ldots\oplus D_r$ is integrable and defines a Riemannian foliation with totally geodesic~leaves;

\noindent\ \
$(b)$ the leaves of $[\xi]\oplus D_0$ are $(2p+1)$-dimensional Sasakian manifolds.
\end{theorem}

\begin{proof}
Consider a unit eigenvector $X$ of $h^2$ with eigenvalue $\mu$. Then, by \eqref{E-c-01} and \eqref{E-nS-01e},
\[
 h^2(\nabla_X\,\xi) = h^2(h-\varphi)X = (h-\varphi)h^2 X = \mu (h-\varphi)X = \mu \nabla_X\,\xi,
\]
hence $\nabla_X\,\xi\in D_\mu$.
On the other hand, \eqref{E-nS-11} and \eqref{E-nabla-xi-Q} imply $\nabla_\xi\,h^2 = 0$, and thus $\nabla_\xi X$ is also an eigenvector of $h^2$ with eigenvalue $\mu$.

Taking unit vectors $X,Y\in D_0$ and applying \eqref{E-nS-11} and $\widetilde Q X=0$, we get $h^2(\nabla_X\,Y) = 0$; thus $\nabla_X\,Y\in [\xi]\oplus D_0$.
Hence, the distribution $[\xi]\oplus D_0$ is integrable with totally geodesic leaves.

If $\mu\ne0$, then taking unit vectors $X,Y\in D_\mu$ and applying \eqref{E-nS-11} and $hY=\mu Y = \varphi Y$, we get
%\begin{equation}\label{E-nS-h2}
% h^2(\nabla_X\,Y) = \mu\,\nabla_X Y - (\nabla_X\,h^2)Y = \mu\,\nabla_X Y - \mu\,g(X, (h-\varphi)Y)\,\xi.
%\end{equation}
%Since $hY=\mu Y = \varphi Y$, from \eqref{E-nS-h2} we get
$h^2(\nabla_X\,Y) = \mu\,\nabla_X Y$.
Hence, each distribution $[\xi]\oplus D_i\ (i = 1,\ldots, r)$ is integrable with totally geodesic leaves.

Let $p>0$.
Using $h^2 \varphi^2=\varphi^2 h^2$, see \eqref{E-nS-01e}, we get for $\mu\ne0$:
\[
 h^2(\varphi^2\nabla_X Y) = \varphi^2(h^2\nabla_X Y) = \mu\,\varphi^2(\nabla_X Y).
\]
Thus $\varphi^2\nabla_X Y \in D_\mu$.
Similarly, using \eqref{E-nS-01d}, we get $\widetilde Q\nabla_X Y \in D_\mu$.
Using the above and \eqref{E-nS-2.1}, we~get
\[
 \nabla_X Y = -\widetilde Q\,\nabla_X Y -\varphi^2\nabla_X Y + \eta(\nabla_X Y)\,\xi,
\]
hence $\nabla_X Y$ belongs to the distribution $[\xi]\oplus D_\mu$ with $\mu\ne0$. By \eqref{E-nS-11} we get
\begin{align*}
 g(\nabla_X Y, Z) & = -(1/\lambda_j^2)\,g(\nabla_X(h^2 Y), Z) = -(1/\lambda_j^2)\,g((\nabla_Xh^2) Y + h^2\nabla_X Y, Z) \\
 & = -(1/\lambda_j^2)\,g(\nabla_X Y, h^2 Z) = 0,\quad
 X\in D_i,\ Y\in D_j,\ Z\in [\xi]\oplus D_0,
\end{align*}
thus the distribution $[\xi]\oplus D_1\oplus\ldots\oplus D_r$ defines a totally geodesic foliation.

First, for any $Z, Z'\in
%[\xi]\oplus
D_0$, we have $(\pounds_\xi\,g)(Z, Z')=0$ since $\xi$ is Killing.
Next, since the distribution $[\xi]\oplus D_0$ is integrable with totally geodesic leaves,
for any
%$Z\in [\xi]\oplus D_0$, $Z'\in D_0$ and
$X\in D_1\oplus\ldots\oplus D_r$ we~get
\[
 g(\nabla_Z X, Z') = -g(\nabla_Z Z', X) = 0.
\]
By the above, we conclude~that for any $Z, Z'\in D_0$ and any $X\in D_1\oplus\ldots\oplus D_r$:
\[
 (\pounds_X\,g)(Z, Z') = g(\nabla_Z X, Z') + g(\nabla_{Z'} X, Z) = 0.
\]
Thus, the distribution $[\xi]\oplus D_1\oplus\ldots\oplus D_r$ defines a Riemannian foliation, hence, a).

If $X$ is an eingenvector of $h^2$ orthogonal to $\xi$ with eigenvalue $\mu$, then also $\varphi X$ is an eingenvector with the same eigenvalue $\mu$.
Hence, the eigenvalue $0$ has odd multiplicity $2p+1$ for some integer $p\ge0$.
If $p>0$,  then the structure $(\varphi,Q,\xi,\eta,g)$ induces a weak nearly Sasakian structure on the totally geodesic
leaves of the distribution $[\xi]\oplus D_0$ whose associated tensor $h$ vanishes.
Thus, $\varphi X = -\nabla_X\,\xi$ for all $X\in [\xi]\oplus D_0$,
and, by Lemma~\ref{L-04}, $QX=X$ for all $X\in[\xi]\oplus D_0$.
Therefore, the induced structure is Sasakian, hence, b).
\end{proof}

\section{Proof of Proposition~\ref{P-R01}}
\label{sec:app}

%\textbf{1}.
The curvature tensor of a weak nearly Sasakian manifold satisfies the equality
\begin{align}\label{E-3.4}
\nonumber
 & g(R_{\,\varphi X,Y}Z, V) +g(R_{\,X,\varphi Y}Z, V) +g(R_{\,X,Y}\,\varphi Z, V) +g(R_{\,X,Y}Z, \varphi V) \\
\nonumber
 & = g(Y,V)\,g((h-\varphi)X,Z)-g(X,Y)\,g(Z, (h-\varphi)V) + g(Y,Z)\,g(X, (h-\varphi)V) \\
 &-(1/2)\,g(Z,V)\,g((h-\varphi)X,Y) +(1/2)\,g(X,Z)\,g(Y, (h-\varphi)V)
\end{align}
for all $X,Y,Z,V\in\mathfrak{X}_M$.
The proof of \eqref{E-3.4} is similar to the proof of equation (19) in \cite{rov-2023b} and we omit it.
%\noindent \textbf{2}.
In the proof of \eqref{E-nS-04ccc}, we also use the following ``small" tensor:
\begin{equation*}
 \delta(X,Y,Z,V) = g(R_{\,X, Y}\widetilde Q Z, V) +g(R_{\,X, Y}Z, \widetilde Q V) - g(R_{\widetilde Q X, Y}Z, V) - g(R_{\,X, \widetilde Q Y}Z, V) .
\end{equation*}
The tensor $\delta$ of a weak nearly Sasakian manifold has some symmetries:
\begin{align*}
 \delta(Y,X, Z,V) = \delta(X,Y, V,Z) = \delta(Z,V, X,Y) = -\delta(X,Y, Z,V) .
\end{align*}
If \eqref{E-nS-04c} is true, then by \eqref{E-3.23b}, we get
\begin{align*}
 & \delta(X,Y,Z,\xi)
 = g(R_{\,X, Y}\,\widetilde Q Z, \xi) -g(R_{\widetilde Q X, Y}Z, \xi) -g(R_{\,X, \widetilde Q Y}Z, \xi) \\
 &\qquad
 = g(R_{\,\xi,\, \widetilde Q Z}Y, X) +g(R_{\,\xi, Z}\,\widetilde Q X, Y) +g(R_{\,\xi, Z}\,X, \widetilde Q Y) = 0 ,
\end{align*}
hence, $\delta(\xi,Y,Z,V)=\delta(X,\xi,Z, V)=\delta(X,Y,\xi, V)=\delta(X,Y,Z, \xi)=0$.

\smallskip

Replacing $X$ by $\varphi X$ in \eqref{E-3.4} and using \eqref{E-nS-2.1}, we have
\begin{align}\label{E-3.13}
\nonumber
 & - g(R_{\,Q X,Y}Z, V) +\eta(X) g(R_{\,\xi,Y}Z, V) +g(R_{\,\varphi X, \varphi Y}Z, V) \\
\nonumber
 & + g(R_{\,\varphi X, Y}\,\varphi Z, V) +g(R_{\,\varphi X, Y}Z, \varphi V) = g(V,Y)\,g((h-\varphi)\varphi X,Z) \\
\nonumber
 & -g(\varphi X,Y)\,g(Z,(h-\varphi)V) + g(Y,Z)\,g(\varphi X, (h-\varphi)V) \\
 &\ - (1/2)\,g(Z,V)\,g((h-\varphi)\varphi X,Y) +(1/2)\,g(\varphi X,Z)\,g(Y, (h-\varphi)V) .
\end{align}
Exchanging $X$ and $Y$ in \eqref{E-3.13}, we find
\begin{align}\label{E-3.14}
\nonumber
 & g(R_{\,X, Q Y}Z, V) +\eta(Y) g(R_{\,\xi,X}Z, V) - g(R_{\,\varphi X, \varphi Y}Z, V) \\
\nonumber
 & + g(R_{\,\varphi Y, X}\,\varphi Z, V) + g(R_{\,\varphi Y, X}Z, \varphi V) = g(V,X)\,g((h-\varphi)\varphi Y,Z) \\
\nonumber
 & - g(\varphi X,Y)\,g(Z, (h-\varphi)V) - g(X,Z)\,g(\varphi(Y, h-\varphi)V) \\
 &\ -(1/2)\,g(Z,V)\,g(X, (h-\varphi)\varphi Y) +(1/2)\,g(Z, \varphi Y)\,g(X, (h-\varphi)V) .
\end{align}
Subtracting \eqref{E-3.14} from \eqref{E-3.13}, we obtain
\begin{align}\label{E-3.15}
\nonumber
 & 2\,g(R_{\,\varphi X,\varphi Y}Z, V) -2\,g(R_{\,X,Y}Z, V) +\eta(X)\,g(R_{\,\xi,Y}Z, V) -\eta(Y)\,g(R_{\,\xi,X}Z, V) \\
\nonumber
 &\ + g(R_{\,\varphi X, Y}\,\varphi Z, V) - g(R_{\,\varphi Y, X}\,\varphi Z, V) + g(R_{\,\varphi X, Y}Z, \varphi V) - g(R_{\,\varphi Y, X}Z, \varphi V) \\
\nonumber
 &\ - g(R_{\widetilde Q X, Y}Z, V) - g(R_{\,X, \widetilde Q Y}Z, V) \\
%%%%%%%%%%%%%
 \nonumber
 & = g(Y,V)\,g((h-\varphi)\varphi X,Z) + g(Y,Z)\,g((h-\varphi)V,\varphi X) + g(Z,V)\,g(\widetilde Q X,Y) \\
 \nonumber
 &\ +(1/2)\,g(\varphi X,Z)\,g(Y, (h-\varphi)V) -(1/2)\,g(\varphi Y,Z)\,g(X, (h-\varphi)V) \\
%%%%%%%%
% \nonumber
 &\ -g(X,V)\,g((h-\varphi)\varphi Y,Z) + g(X,Z)\,g(Y, \varphi(h-\varphi)V) .
\end{align}
Then, replacing $Z$ by $\varphi Z$ and also $V$ by $\varphi V$ in \eqref{E-3.4} and using \eqref{E-nS-2.1}, we get two equations
\begin{align}\label{E-3.16}
\nonumber
 & - g(R_{\,X, Y}Q Z, V) = -\eta(Z)\,g(R_{\,X,Y}\,\xi, V) - g(R_{\,X,Y}\,\varphi Z, \varphi V) - g(R_{\,X,\varphi Y}\,\varphi Z, V) \\
%%%%
\nonumber
 &\ - g(R_{\,\varphi X,Y}\,\varphi Z, V) - g(Y,V)\,g(\varphi(h-\varphi)X, Z)-g(X,Y)\,g(\varphi Z, (h-\varphi)V) \\
\nonumber
 &\  + g(Y,\varphi Z)\,g(X, (h-\varphi)V) -(1/2)\,g(\varphi Z,V)\,g((h-\varphi)X,Y) \\
 &\ +(1/2)\,g(X,\varphi Z)\,g(Y, (h-\varphi)V),\\
\label{E-3.17}
\nonumber
 & -g(R_{\,X, Y}Z, Q V) = -\eta(V)\,g(R_{\,X,Y}Z,\,\xi) - g(R_{\,X,Y}\,\varphi Z, \varphi V) - g(R_{\,\varphi X, Y}Z, \varphi V) \\
%%%%
 \nonumber
 &\ - g(R_{\,X,\varphi Y}Z, \varphi V)  + g(Y,\varphi V)\,g((h-\varphi)X,Z)-g(X,Y)\,g(Z, (h-\varphi)\varphi V) \\
 \nonumber
 &\  + g(Y,Z)\,g(X, (h-\varphi)\varphi V) -(1/2)\,g(Z,\varphi V)\,g((h-\varphi)X,Y) \\
 &\ +(1/2)\,g(X,Z)\,g(Y, (h-\varphi)\varphi V) .
\end{align}
Adding \eqref{E-3.16} to \eqref{E-3.17}, we get
\begin{align*}
 &-2\,g(R_{\,X, Y}Z, V) = -2\,g(R_{\,X,Y}\,\varphi Z, \varphi V) + g(R_{\,X, Y}\widetilde Q Z, V) + g(R_{\,X, Y}Z, \widetilde Q V) \\
 &\ -\,\eta(Z)\,g(R_{\,X,Y}\,\xi, V)
   - g(R_{\,X,\varphi Y}\,\varphi Z, V) - g(R_{\,\varphi X,Y}\,\varphi Z, V)\\
 &\ -\,\eta(V)\,g(R_{\,X,Y}Z,\,\xi) - g(R_{\,\varphi X, Y}Z, \varphi V) - g(R_{\,X,\varphi Y}Z, \varphi V) \\
%%%%%
\nonumber
 &\ - g(Y,V)\,g(\varphi(h-\varphi)X, Z) + g(Y,\varphi Z)\,g(X, (h-\varphi)V) \\
%%%%
\nonumber
 &\ + g(Y,\varphi V)\,g((h-\varphi)X,Z) +2\,g(X,Y)\,g(Z, \widetilde Q V) + g(Y,Z)\,g(X, (h-\varphi)\varphi V) \\
 &\ +(1/2)\,g(X,\varphi Z)\,g(Y, (h-\varphi)V) +(1/2)\,g(X,Z)\,g(Y, (h-\varphi)\varphi V) .
\end{align*}
Substituting the above equation into \eqref{E-3.15}, we have
\begin{align}\label{E-3.18}
\nonumber
 & 2\,g(R_{\,\varphi X,\varphi Y}Z, V) -2\,g(R_{\,X,Y}\,\varphi Z, \varphi V) - \eta(Z)\,g(R_{\,\xi, V}X,Y) + \eta(V)\,g(R_{\,\xi,Z}\,X,Y) \\
\nonumber
 & +\,\eta(X)\,g(R_{\,\xi,Y}Z, V) -\eta(Y)\,g(R_{\,\xi,X}Z, V) +\delta(X,Y,Z,V) \\
%%%
\nonumber
 & - g(Y,V)\,g(\varphi(h-\varphi)X, Z) + g(Y,\varphi Z)\,g(X, (h-\varphi)V) \\
%%%%
\nonumber
 & + g(Y,\varphi V)\,g((h-\varphi)X,Z) + g(Y,Z)\,g(X, (h-\varphi)\varphi V) + 2\,g(X,Y)\,g(Z, \widetilde Q V) \\
 & + (1/2)\,g(X,\varphi Z)\,g(Y, (h-\varphi)V) +(1/2)\,g(X,Z)\,g(Y, (h-\varphi)\varphi V) = 0 .
\end{align}
 Replacing $X$ by $\varphi X$ and also $Y$ by $\varphi Y$ in \eqref{E-3.18} and using \eqref{E-nS-2.1}, we obtain
\begin{align}\label{E-3.19}
\nonumber
 &2\,g(R_{\,Q X, Q Y}Z, V) - 2\,g(R_{\,\varphi X,\varphi Y}\,\varphi Z,\varphi V) - 2\,\eta(X)\,g(R_{\,\xi, Q Y}Z, V) \\
\nonumber
 & + 2\,\eta(Y)\,g(R_{\,\xi,Q X}Z, V)
  - \eta(Z)\,g(R_{\,\xi, V}X,Y) + \eta(V)\,g(R_{\,\xi, Z}\,\varphi X, \varphi Y) \\
%%%%%%%%%%%
\nonumber
& + \delta(\varphi X, \varphi Y,Z,V)
 + g(\varphi Y,V)\,g((h-\varphi)\varphi X,\varphi Z) - g(\varphi Y,\varphi Z)\,g(\varphi(h-\varphi)V, X) \\
\nonumber
 & +g(\varphi Y,\varphi V)\,g((h{-}\varphi)\varphi X,Z)
 + g(\varphi Y,Z)\,g(\varphi X, (h{-}\varphi)\varphi V) -2\,g(\varphi^2 X, Y)\,g(Z, \widetilde Q V)\\
%\nonumber
 &  +(1/2)\,g(\varphi^2 X, Z)\,g(Y, \varphi(h-\varphi)V) +(1/2)\,g(\varphi X,Z)\,g(\varphi Y, (h-\varphi)\varphi V)  =0 .
\end{align}
Replacing $V$ by $\xi$ in \eqref{E-3.19}, we obtain
\begin{align}\label{E-3.20}
\nonumber
 & 2\,g(R_{\,Q X, Q Y}Z, \xi) -2\,\eta(X)\,g(R_{\,\xi, Q Y}Z, \xi) +2\,\eta(Y)\,g(R_{\,\xi,Q X}Z, \xi)  \\
 & +\,g(R_{\,\xi, Z}\,\varphi X, \varphi Y)  +\delta(\varphi X, \varphi Y,Z, \,\xi) = 0.
\end{align}
 Replacing $X$ by $\varphi X$ and also $Y$ by $\varphi Y$ in \eqref{E-3.20}, we obtain
\begin{align}\label{E-3.21}
\nonumber
 & 4\,g(R_{\,Q\,\varphi X,\, Q\,\varphi Y}Z, \xi) + 2\,\delta(\varphi^2 X, \varphi^2 Y, Z, \,\xi) +2\,g(R_{\,\xi, Z}\,Q X, Q Y) \\
 & -2\,\eta(X)\,g(R_{\,\xi, Z}\,\xi, Q Y) +2\,\eta(Y)\,g(R_{\,\xi, Z}\,\xi, Q X) = 0.
\end{align}
Adding \eqref{E-3.20} and \eqref{E-3.21} and applying $Q = {\rm id}_{\,TM} + \widetilde Q$, we get
\begin{align}\label{E-3.22}
\nonumber
 3\,g(R_{\,\xi, Z}\,\varphi X, \varphi Y) = -4\,g(R_{\,\xi, Z}\,\widetilde Q\varphi X, \varphi Y) -4\,g(R_{\,\xi, Z}\,\varphi X, \widetilde Q\varphi Y) \\
 -4\,g(R_{\,\xi, Z}\,\widetilde Q\varphi X, \widetilde Q\varphi Y)
 + \delta(\varphi X, \varphi Y, Z, \,\xi) +2\,\delta(\varphi^2 X, \varphi^2 Y, Z, \,\xi).
\end{align}
Using the condition \eqref{E-nS-04c}, we get
\[
 g(R_{\,\xi, Z}\,\widetilde Q\varphi X, \varphi Y) = g(R_{\,\xi, Z}\,\varphi X, \widetilde Q\varphi Y)
 =g(R_{\,\xi, Z}\,\widetilde Q\varphi X, \widetilde Q\varphi Y) =0
\]
and $\delta(\varphi X, \varphi Y, Z, \,\xi) = \delta(\varphi^2 X, \varphi^2 Y, Z, \,\xi)=0$.
Therefore, from \eqref{E-3.22} we obtain \eqref{E-nS-04ccc}.
Since the restriction $\varphi\,|_{\,\ker\eta}$ is non-degenerate, the distribution $\ker\eta$ is curvature invariant, see \eqref{E-nS-04cc}.
\hfill$\square$

\section{Conclusions}

We have shown that the weak nearly Sasakian structure is useful for studying contact metric structures and totally geodesic foliations.
Some well known results on nearly Sasakian manifolds were extended
to weak nearly Sasakian manifolds with the conditions \eqref{E-nS-10} and \eqref{E-nS-04c},
and our main result is that such weak nearly Sasakian manifolds have a foliated structure with two types of totally geodesic foliations.
Motivated by a criterion for an almost contact metric manifold to be Sasakian in \cite{NDY-2018}, we pose the following question:
{are weak nearly Sasakian manifolds of dimension higher than five, satisfying the conditions \eqref{E-nS-10} and \eqref{E-nS-04c}, Sasakian manifolds}?
Based on the numerous applications of the nearly Sasakian structure, we expect that certain weak structure will also be useful for geometry and physics,
e.g., in twistor string theory.

\end{document}